\documentclass[11pt,a4paper]{amsart}
\usepackage[english]{babel}
\usepackage{pst-grad} 
\usepackage{pst-plot} 
\usepackage{pstricks}
\usepackage{amsmath,amssymb,mathrsfs,amsthm}
\usepackage[dvips]{graphicx}
\usepackage{epsfig}

\newcommand{\RR}{\mathbb R}

\newcommand{\TT}{\mathbb T}
\newcommand{\pat}{\partial_t}
\newcommand{\pax}{\partial_x}

\newcommand{\jeps}{\mathcal{J}_\varepsilon*}

\newcommand{\vertiii}[1]{{\left\vert\kern-0.25ex\left\vert\kern-0.25ex\left\vert #1 
    \right\vert\kern-0.25ex\right\vert\kern-0.25ex\right\vert}}

\usepackage[pdfauthor={Rafael Granero-Belinchon},pdftitle={...},bookmarks,colorlinks]{hyperref}

\newcounter{comentcount}
\newcounter{teocount}
\newtheorem{ass}{Assumption}

\newtheorem{corol}{Corollary}
\newtheorem{teo}[teocount]{Theorem}  
\newtheorem{defi}{Definition}

\newtheorem{rem}{Remark}

\title[Chemotaxis model with nonlocal and semilinear flux]{Boundedness of large-time solutions to a chemotaxis model with nonlocal and semilinear flux.}

\author[J. Burczak]{Jan Burczak}
\email{jb@impan.pl}
\address{Institute of Mathematics of the Polish Academy of Sciences, Warsaw, 21 00-956, Poland}

\author[R. Granero-Belinch\'{o}n]{Rafael Granero-Belinch\'{o}n}
\email{rgranero@math.ucdavis.edu}
\address{Department of Mathematics, University of California, Davis, CA 95616, USA}

\begin{document}

\begin{abstract}
A semilinear version of parabolic-elliptic Keller-Segel system with the \emph{critical} nonlocal diffusion is considered in one space dimension. We show boundedness of weak solutions under very general conditions on our semilinearity. It can degenerate, but has to provide a stronger dissipation for large values of a solution than in the critical linear case or we need to assume certain (explicit) data smallness. Moreover, when one considers a logistic term with a parameter $r$, we obtain our results even for diffusions slightly weaker than the critical linear one and for arbitrarily large initial datum, provided  $r >1$. For a mild logistic dampening, we can improve the smallness condition on the initial datum up to $\sim \frac{1}{1-r}$.
\end{abstract}

\maketitle 


\section{Introduction}
In this paper we study the following model
\begin{eqnarray}\label{eqa1}
\pat u & =&  \pax(-\mu(u)Hu+u\pax v) +ru(1-u),\; \;x\in\TT, \, t\in\RR^+ \\ 
\label{eqa2}
\pax^2 v & =&  u-\langle u \rangle,\; \;x\in\TT, \, t\in\RR^+,
\end{eqnarray}
where $u=u(x,t), v=(x,t)$, $H$ stands for the (periodic) Hilbert transform, \emph{i.e.}
$$
\widehat{Hu}(\xi)=-i\frac{\xi}{|\xi|}\hat{u}(\xi).
$$
$\TT=[-\pi,\pi]$, $r \ge 0$ and $\mu$ is a certain function (semilinearity), precised in what follows. Before formulating our results let us explain our motivations to study the system \eqref{eqa1}-\eqref{eqa2}.
\subsection{Motivation} (a -- \emph{mathematical biology}) One of the basic systems studied in the context of chemotaxis is the parabolic-elliptic Keller-Segel system (also known as the Smoluchowski-Poisson system)
\begin{equation}\label{eqa1local}
\pat u  =  \nabla\cdot(\mu\nabla u-u\nabla \phi),  \quad x\in\TT, \, t\in\RR^+,
\end{equation}
where $d\geq1$ denotes the spatial dimension, $\TT^d=[-\pi,\pi]^d$, $\mu>0$ is a constant and $\phi$ is recovered from $u$ through some operator, \emph{i.e.} $\phi(x,t)=T(u(x,t))$. 
In many cases $\phi$ satisfies the Poisson equation
\begin{equation}\label{eqa2local}
- \Delta \phi  =  u-\langle u \rangle, \quad \;x\in\TT^d, t\in\RR^+.
\end{equation}

In this notation, $u$ represents the concentration of cells, $\langle u \rangle$ its space average and $\phi$ gives us the concentration of a chemical substance that attracts cells. It is biologically justified to enrich the equation  \eqref{eqa1local} with the logistic term, obtaining
\begin{equation}\label{eqa1localR}
\pat u  =  \nabla\cdot(\mu\nabla u-u\nabla \phi) +ru(1-u),\;x\in\TT^d, t\in\RR^+,
\end{equation}
where $r \ge 0$. The model \eqref{eqa1localR}-\eqref{eqa2local} is related to the parabolic-elliptic simplification of  the cell kinetics model M$8$ in \cite{Hillen3}, that describes a  bacterial pattern formation or  cell movement and growth during angiogenesis.
 
 Another application of the model \eqref{eqa1localR}-\eqref{eqa2local} occurs in tumor growth. In particular, this model is related to the three-component urokinase plasminogen invasion model (see \cite{Hillen1}). There is a huge literature on the mathematical study of a numerous versions of \eqref{eqa1localR}-\eqref{eqa2local} in the context of mathematical biology, see \cite{Bi6, blanchet2010functional, BCM, Dolbeault4, Carrillo2, Dolbeault2, jager1992explosions} and the references therein.\\
 (b -- \emph{natural sciences}) Let us take in  \eqref{eqa1localR}-\eqref{eqa2local} $v := - \phi$. The resulting system 
 \begin{equation}\label{eqa1local'}
\pat u  =  \nabla\cdot(\mu\nabla u+u\nabla v) +ru(1-u),\;x\in\TT^d, t\in\RR^+,
\end{equation}
 \begin{equation}\label{eqa2local'}
 \Delta v  =  u-\langle u \rangle,\;x\in\TT^d, t\in\RR^+
\end{equation}
in the case $r=0$ is important in mathematical cosmology and gravitation theory. It is very similar in spirit to the Zel'dovich approximation used in cosmology to study the formation of large-scale structure in the primordial universe, see also \cite{AGM, Bi1}. It is also connected with the Chandrasekhar equation for the gravitational equilibrium of polytropic stars, statistical mechanics and the Debye system for electrolytes, see \cite{BilNad94}.
\begin{rem}
In what follows, we consider a system with the sign $"+"$ in front of terms $u\nabla v$ and $\Delta v$, compare motivation (b), remembering that letting $\phi :=-v$ we get the equations studied in mathematical biology (see motivation  (a)).
\end{rem}
\subsection{Central problem}
The focal point in the studies of  solutions to \eqref{eqa1local'}-\eqref{eqa2local'} is the matter of distinguishing between the blowup and global-in-time regimes in correlation with the dimension $d$, initial data and parameters of the system. Roughly speaking, it turns out that with $\mu =1$, $r=0$ there is a $8 \pi$ criticality of the initial mass $\|u_0 \|_{L^1}$ in two dimensions. Below this threshold one can have the global existence of bounded solutions and above it there is a finite time blowup  (see for instance \cite{BouCalChapter}). In the one-dimensional case, the diffusion  $\nabla\cdot(\nabla u)$ is strong enough to give the global existence. On the other hand, for $d >2$ it is too weak (let us remark here that the logistic term $ru (1-u)$ generally helps the global existence, compare \cite{WinklerLogistic}). In this context it is mathematically interesting to find, for a fixed dimension $d$, a "critical" diffusive operator that sits on the borderline of the blowup and global-in-time regimes. There are at least two approaches to this problem, both  justified also from the point of view of applications. One is to consider semilinear diffusion $\nabla \cdot (\mu ( u) \nabla u)$, see for instance \cite{Bertozzi, blanchet2009critical,  CieslakStinner, BurczakCieslak, CieslakLaurencot, TaoWinkler}. Another one is to replace the standard diffusion with the fractional one. We focus on this case. Let us consider 
\begin{equation}\label{eqa1nolocal}
\pat u  =  -\mu\Lambda^\alpha u+\nabla\cdot(u\nabla v) +ru(1-u),\quad x\in\TT^d, t\in\RR^+,
\end{equation}
 \begin{equation}\label{eqa2nolocal}
\Delta v  =  u-\langle u \rangle,\quad x\in\TT^d, t\in\RR^+,
\end{equation}
where $\mu>0$ is a constant and the operator $\Lambda^{\alpha}$ is defined using the Fourier transform
$$
\widehat{\Lambda^\alpha u}(\xi)=|\xi|^\alpha\hat{u}(\xi).
$$
It turns out that with $r=0$ and $d=1$ there are global-in-time solutions for $\alpha>1$ and blowups for $\alpha<1$, compare \cite{bournaveas2010one} and \cite{escudero2006fractional} (see also \cite{AGM}). The case $\alpha=1$ seems critical and, to the best of our knowledge, the sharpest results up to now is the global boundedness for small data. In particular, it is shown in \cite{bournaveas2010one} that there exists a constant $K$, such that $\|u_0\|_{L^1}\leq K$ implies the global existence of solutions. Later on, in \cite{AGM}, the authors proved that $\|u_0\|_{L^1}\leq 1/(2\pi)$ implies global existence and the convergence towards the homogeneous steady state. In this context we refer also to  \cite{BurczakGranero, li2010exploding}. 

In this paper we propose a slight semilinear strengthening of the diffusion in \eqref{eqa1nolocal} that provides bounded solutions for any $r \ge 0$. We also study the regularization due the logistic dampening for the diffusions equal to that of \eqref{eqa1nolocal} or slightly weaker, compare example \eqref{example}. 

The case of  \eqref{eqa1}-\eqref{eqa2} with $\alpha=1$, $r = 0$ and $\mu (s)= s + \nu$ appears in \cite{GO}, where the authors address the local/global existence and the qualitative behaviour of the solutions. Similar equations have been studied in \cite{BaeGranero, Carrillo, CC, CC2, CCCF} in the context of fluid dynamics. In particular, the equation
\begin{equation}\label{eq2angel}
\pat u=\pax\left(-[u+\nu] Hu\right).
\end{equation}
has been proposed as a one-dimensional model of the 2D Vortex Sheet problem or the 2D surface quasi-geostrophic equation. Notice that \eqref{eqa1} reduces to \eqref{eq2angel} when $v\equiv0,r=0$ and $\mu(x)=x+\nu$.
\subsection{Basic notation and plan of the paper}

We write $H^s(\TT)$ for the usual $L^2$-based Sobolev spaces with norm
$$
\|f\|_{H^s}^2:=\|f\|_{L^2}^2+\|f\|_{\dot{H}^s}^2, \|f\|_{\dot{H}^s}:=\|\Lambda^s f\|_{L^2},
$$
and
$$
\langle u\rangle:=\frac{1}{2\pi}\int_{\TT}u(x)dx.
$$
For a given initial data $u_0$, we introduce the following notation
\begin{equation}\label{n1}
\mathcal{N}_1:=\max\{2\pi,\|u_0\|_{L^1}\}
\end{equation}

Notice that for periodic functions, the half-laplacian in one dimension has the following kernel representation
$$
\Lambda f(x)=\frac{\text{p.v.}}{2\pi}\int_\TT\frac{f(x)-f(y)}{\sin^2\left(\frac{x-y}{2}\right)}dy.
$$
The remainder of this paper is organized as follows: In Section \ref{S1} we present the statement of our results. In Section \ref{S2} we prove Theorem \ref{teo1}. In Section \ref{sec4} we prove Theorem \ref{teoRlarge}. Finally, in Section \ref{S3} we present our proof of Theorem \ref{teo3}.
\section{Statement of results}\label{S1}
Given  the initial data $u_0(x)\geq0$, we have the following definition of a~weak solution to the system \eqref{eqa1}-\eqref{eqa2}:

\begin{defi}\label{defipe}
Choose $u_0 \in L^2(\TT)$. Fix arbitrary $T \in (0, \infty)$. The couple \[(u,v)\in L^\infty(0,T;L^2(\TT))\times L^\infty(0,T;H^{1}(\TT))\] is a solution of \eqref{eqa1}-\eqref{eqa2} if and only if
$$
\int_0^T\int_\TT -\pat\phi \, u+\pax\phi \, (-\mu(u)Hu+u\pax v) -\phi \, ru(1-u)dxdt - \int_\TT \phi(x,0) u_0dx=0,
$$
$$
\int_0^T\int_\TT \pax \varphi \,  \pax v + \varphi \, (u-\langle u\rangle)dxdt=0,
$$
for every test function $\phi(x,t),\varphi(x,t)\in C^\infty((-1,T)\times\TT)$ with a compact support in time and periodic in space.

\end{defi}
\begin{defi}\label{defipeT}
If a solution $(u,v)$ verifies Definition \ref{defipe} for any $T < \infty$, we call it a \emph{large-time  weak solution}.
\end{defi}

We will use the following entropy (or free energy) functional
\begin{equation}\label{entropy}
\mathcal{F}(u(t))=\int_{\TT} u(t)\log(u(t))-u(t)+1.
\end{equation}

\subsection{Case of linearly degenerating $\mu$ and any $r \ge 0$}
The results presented in this section do not use extra information in estimates that follows from the logistic term $r u (1-u)$. Hence they hold for any value of  $r \ge 0$, including small ones. 
The semilinearity $\mu$ of our diffusion will be generated here by a function $\gamma$ as follows
\begin{equation}\label{Gamma}
 \mu(s) := \gamma(s) s.
\end{equation}
Let us introduce also
\[
 \Gamma(s):= \int_{0}^s\gamma(y)dy.
\]
We work within
\begin{ass}\label{ass1} The semilinearity $\mu$ is differentiable and its derivative $\mu'$  is bounded for bounded arguments, i.e. there exists a finite function $C$ such that
\[
\mu' (s) \le C(s)
\]
for any $s \in [0, \infty)$. Moreover, $\gamma$ of \eqref{Gamma} satisfies for any $y \in [0, \infty)$
\begin{equation}\label{mu}
\gamma(y) \ge \delta>0
\end{equation}
for a fixed $\delta > 0$ and there exists $0\leq y_0<\infty$ such that 
\begin{equation}\label{assA}
 \gamma(y)\geq1\text{ for }y\geq y_0. 
 \end{equation}
\end{ass}
The fact that $\mu$ is \emph{linearly degenerating} is understood in the sense of condition \eqref{mu}, as it allows for $\mu (s) = \delta s$ for small $s$.
\begin{teo}\label{teo1}Let $0\leq u_0\in L^\infty$ be the initial data for \eqref{eqa1}-\eqref{eqa2} under Assumption \ref{ass1}. Then there exists at least one global in time weak solution to  \eqref{eqa1}-\eqref{eqa2} (in the sense of Definitions \ref{defipe},  \ref{defipeT}). Furthermore, this solution enjoys additionally the following regularity
$$
u \in L^\infty(0,T; L^\infty(\TT))\cap L^2(0,T;H^{0.5}(\TT)) \quad \forall \; T< \infty,
$$
where the $L^\infty$ bound is $T$-independent.
\end{teo}
\subsection{Results using the logistic dampening}
Now we formulate a result that allows for the the critical linear nonlocal diffusion (\emph{i.e.} $\mu \equiv c$) at the cost of using a relation between the lower bound on $\mu$, the initial mass $\langle u_0\rangle$ and $r$. Moreover, for strictly positive data it  generalizes Theorem \ref{teo1} over any $\mu (s)$ that is positive for $s>0$. Particularly, we do not need to assume here the sublinear profile of degeneration of $\mu$.

Now we assume the following hypothesis on the semilinearity.
\begin{ass}\label{ass2} 
The semilinearity $\mu$ is differentiable and its derivative $\mu'$  is bounded for bounded arguments, i.e. there exists a finite function $C$ such that
\[
\mu' (s) \le C(s)
\]
for any $s \in [0, \infty)$. Moreover $\mu$ is positive for positive arguments, \emph{i.e.}
$$
\mu(x_0)=0\Rightarrow x_0=0,
$$
and there exist $\delta \ge 0$ such that 
\[
\mu (s) \ge \delta.
\]
\end{ass}
Observe that above  we allow $\delta = 0$ (then the condition of positivity for positive arguments prevails).

In the following  result  $  \delta \ge 0$ comes from Assumption \ref{ass2}  and $r$ is the parameter of the logistic term.
\begin{teo}\label{teoRlarge}
Let $0\leq u_0\in L^\infty$. If, in addition to Assumption \ref{ass2} we have
\begin{equation}\label{c:teoL}
r +  \delta (4\pi^2\max\{\langle u_0\rangle,1\})^{-1} > 1, 
\end{equation}
and either 
\begin{equation}\label{mularger}
\delta > 0
\end{equation}
or 
\begin{equation}\label{mularger'}
 \text{ess}\,\text{min}_x u_0 >0
\end{equation}
then there exists at least one global in time weak solution to  \eqref{eqa1}-\eqref{eqa2} (in the sense of Definitions \ref{defipe},  \ref{defipeT}). This solution enjoys additionally the following regularity
$$
u \in L^\infty(0,T; L^\infty(\TT))\cap L^2(0,T;H^{0.5}(\TT)) \quad \forall \; T< \infty,
$$
where the $L^\infty$ bound is $T$-independent.
\end{teo}
In the case $r=0$ we can provide a simpler condition, according to
\begin{corol}\label{c1}
In the case $r=0$ Theorem \ref{teoRlarge} is valid with condition \eqref{c:teoL} replaced with
\begin{equation}\label{c:teoL'}
  \frac{ \delta }{4\pi^2 \langle u_0\rangle} > 1. 
\end{equation}
\end{corol}
\subsection{Remarks} The weakest semilinear diffusion allowed by Theorem  \ref{teo1} is $\mu (s) = \delta (s) s$ with: $\delta (s) \ge 1$ for large $s$ and being  arbitrary positive constant otherwise.  This is the previously mentioned {\emph{semilinear strengthening}} of the critical nonlocal linear diffusion $-\mu \Lambda u$ with constant $\mu$. 

On the other hand, in Theorem  \ref{teoRlarge} this critical diffusion is admissible with any $\mu \equiv c$ for arbitrary $c>0$, provided the logistic parameter $r \ge 1$. What's more, for $r>1$ and strictly positive data we can allow for
\begin{equation}\label{example}
\mu (s)  \begin{cases} = 0 &\text{ for } s = 0, \\
> 0 &\text{ for } s > 0, \\
\end{cases}
\end{equation}
hence we do not assume any precise profile of degeneracy at $0$ of our semilinearity $\mu$, hence it can be slightly weaker than the linear $\mu \equiv c>0$.
From the proof of Theorem  \ref{teoRlarge} one sees that even a diffusion that vanishes outside a certain interval of arguments is allowed, compare Remark \ref{rem:stupido}.
Finally, for $r \in [0, 1)$ we need in Theorem  \ref{teoRlarge}  to mitigate the weaker logistic dampening with larger diffusion, namely such that
\[
 \delta  > (1-r) (4\pi^2\max\{\langle u_0\rangle,1\}).
\]
Observe that the bigger the initial mean value, the stronger diffusions we need. See also Corollary  \ref{c1}.

Let us now compare Theorems \ref{teo1},  \ref{teoRlarge} with known results. 
\begin{itemize}
\item To get the system considered  in \cite{GO} we take $r=0$ and $\gamma (x)=1 +  \frac{ \nu }{x}$. This falls under our assumption \eqref{assA}. We see now that our Theorem  \ref{teo1} recovers the result of Theorem 5.2 in \cite{GO} and sharpens it with respect to the admissible initial data. Namely,  we have removed the $H^{0.5}$ smoothness requirement and, more importantly, the smallness assumption $ \|u_0\|_{L^1} \le \frac{2}{3} \nu$ of \cite{GO}.
\item We allow in Theorem  \ref{teo1} for much more general semilinear diffusions $\mu$  than in \cite{GO} and for the logistic term.
\item Theorem   \ref{teoRlarge} is the semilinear version of Theorem 3 in \cite{BurczakGranero}. Additionally, it allows to weight out the admissible $r \ge 0$ with dissipation. Namely, the condition $r \ge \frac{1}{\nu} $ of Theorem 3 in \cite{BurczakGranero} reads in our setting $r \ge 1$, which is exactly the value of $r$ in \eqref{c:teoL} that allows for any lower bound $\delta >0$ for $\mu$. Similarly, the condition  $\delta>2\pi (\max\{|u_0|_{L^1}, 2 \pi \})$ of Theorem 3 in \cite{BurczakGranero} complies with  \eqref{c:teoL} for $r=0$.
\item When restricted to linear case, \emph{i.e.} $\mu \equiv c$, Corollary \ref{c1} is inline with results in \cite{AGM, bournaveas2010one}. Moreover, our condition \eqref{c:teoL'} is explicit and says that the threshold mass for a (debatable) blowup is at least $  \frac{ \delta }{2\pi} $.

\end{itemize}
\subsection{Results for a nonlocal porous medium  type equation}
If we take $v\equiv 0$ and $r=0$ in \eqref{eqa1} we get
\begin{equation}\label{eqa3}
\pat u = \pax(-\mu(u)Hu),\;x\in\TT, t\in\RR^+.
\end{equation}
We have the following result
\begin{teo}\label{teo2}Let $0\leq u_0\in L^\infty$ be the initial data for \eqref{eqa3} with $\mu(s)$ following Assumption \ref{ass1}. Then there exists at least one global in time weak solution $u$ of \eqref{eqa3}. Furthermore, this solution is 
$$
u \in  L^\infty(0,T; L^\infty(\TT))\cap L^2(0,T;H^{0.5}(\TT)) \quad \forall \; T< \infty
$$
\end{teo}
The proof of this theorem is similar to the proof of Theorem \ref{teo1}, so we omit it. 

Finally, we provide also the following result on asymptotics of solutions to a special case of \eqref{eqa3}.
\begin{teo}\label{teo3}Let $u_0\in L^\infty$ such that $\text{ess}\,\text{min}_x u_0 >0$ be the initial data for \eqref{eqa3} with $\mu(x)=x$ and assume that $\langle u_0\rangle=1$. Then the global in time solution $u(x,t)$ of \eqref{eqa3} tends to the homogeneous steady state $u_\infty\equiv 1$ and satisfies
$$
\mathcal{F}(u(t))\leq C(u_0)e^{-2 (\text{ess}\,\text{inf}_x u_0) \,t}.
$$
\end{teo}

In particular this Theorem covers the case $\lambda=0,$ $s=1/2$ in \cite{CarrilloVazquez} for periodic and positive initial data.

\section{Proof of Theorem \ref{teo1}}\label{S2}
\subsection{The approximate problems}
Let's consider a family of Friedrichs mollifiers  $\mathcal{J}_\epsilon$. We define the regularized initial data
$$
u^\epsilon(x,0)=u_0^\epsilon(x)=\jeps u_0(x)\geq0
$$
and consider the approximate problems
\begin{eqnarray}\label{eqa1approx}
\pat u^\epsilon & = & \pax(-\mu(u^\epsilon)Hu^\epsilon+u^\epsilon\pax v^\epsilon) +ru^\epsilon(1-u^\epsilon)+\epsilon\pax^2 u^\epsilon, \\ 
\label{eqa2approx}
\pax^2 v^\epsilon & = & u^\epsilon-\langle u^\epsilon\rangle.
\end{eqnarray}
Each of these problems has a local-in-time smooth solution with the maximal time of existence $T_\epsilon$. This can be shown via a fixed point argument (basically, Picard's Theorem in Banach spaces). For the time being we will work within this local time of existence.
\subsection{$L^1$ bound}
Integrating the equation \eqref{eqa1approx} and using Jensen's inequality as in \cite{BurczakGranero}, we have
\begin{equation}\label{L1}
\|u^{\epsilon}(t)\|_{L^1}\leq \mathcal{N}_1=\max\{\|u_0\|_{L^1},2\pi\}.
\end{equation}

\subsection{Pointwise bounds}
Let's denote the point where $\min_ x u^\epsilon(t)$ is attained as $\underline{x}_t$. Similarly, we write $\overline{x}_t$ for the point where the maximum is reached. In other words,
\begin{equation}\label{notationmax}
\min_{x\in\TT} u^\epsilon(t)=u^\epsilon(\underline{x}_t,t),\;\max_{x\in\TT} u^\epsilon(t)=u^\epsilon(\overline{x}_t,t).
\end{equation}
Then, using the same arguments as in \cite{AGM, BaeGranero, BurczakGranero, GO}, we prove
\begin{equation}\label{min}
\min_{x\in\TT} u^\epsilon(t) \geq  \min_{x\in\TT} u^\epsilon_0 e^{\int_0^t -\gamma(u^{\epsilon}(\underline{x}_s))\Lambda u^{\epsilon}(\underline{x}_s)+u^\epsilon(\underline{x}_s)-\langle u^\epsilon (s)\rangle+r(1-u^{\epsilon}(\underline{x}_s))ds}\geq0
\end{equation}
and, using
$$
\Lambda u^\epsilon(\overline{x}_t,t)\geq u^\epsilon(\overline{x}_t,t)-\langle u^\epsilon(t)\rangle ,
$$
we get the following ODI for $X (t) := \|u^\epsilon(t)\|_{L^\infty}$
\[
\dot X \le X\big((X-\langle u^\epsilon\rangle) (1- \gamma (X))+r(1-X) \big).
\]
Let's assume first that $r>0$. As $X(t)\geq \langle u^\epsilon(t)\rangle$, the sign of $\dot X$ is negative provided $\max\{r(1-X),1-\gamma(X)\}<0$. Recalling Assumption \ref{ass1} (in particular condition \eqref{assA}), we can ensure the existence of $s_0\in\RR^+$ such that $\gamma(s)\geq 1$ if $s\geq s_0$. We can always choose $s_0$ such that $s_0>1$ and $X(0) \le u_0 <s_0$. Then, we have the alternative
   \begin{itemize}
\item[(i)] Either $X (t)  \le s_0$ for all times.
\item[(ii)] Or there exists $t_0 > 0$ such  $X (t_0)  = s_0$ and $X$ crosses $s_0$ for the first time at $t_0$. 
\end{itemize}
In view of the above choices of $s_0$ we see that $1- \gamma (X (t_0)) \le 0$ and $r(1-X(t_0))<0$. Hence our ODI gives $\dot X(t)|_{t=t_0} < 0$, which rules out case (ii). As a consequence, we have $X(t)\leq s_0$ for every time.

In the case $r=0$, we have $\langle u^\epsilon (t) \rangle=\langle u_0 \rangle$. 
As before, we can choose $s_0>\max (1,X(0))$ and such that $1-\gamma(s_0)\leq0$. We have the same alternative (i)-(ii). In the case (i) we have the desired bound. To rule out the case (ii), notice that 
 applying Gronwall's inequality to our ODI we get
\begin{equation}\label{gronwallX}
X(t) \le X(t_0)e^{\int_{t_0}^t(X(\tau)-\langle u_0\rangle) (1- \gamma (X(\tau)))d \tau}\quad\text{for } t\geq t_0.
\end{equation}
According to (ii), $X$ crosses $s_0$ for the first time at $t_0$, so there exists $\delta  > 0$ such that $1-\gamma(X(\tau))\leq0$ for $\tau \in [t_0, t_0 + \delta]$. But then \eqref{gronwallX} gives $X(t) \le s_0$ for $t \in [t_0, t_0 + \delta]$, because $X(t_0) = s_0$ and the exponential function has nonpositive exponent. Hence we have falsified (ii).

Consequently we have the bound
\begin{equation}\label{max}
\|u^\epsilon(t)\|_{L^\infty}\leq s_0 ( \|u_0\|_{L^\infty}, r, \gamma).
\end{equation}

\subsection{$\dot{H}^1$ bound}
For the time being, we have worked within the local time of existence $T_\epsilon$. Now we fix any $T \in (0, \infty)$ and prove that  $T_\epsilon \ge T$ for any $\epsilon$. We test \eqref{eqa1approx} against $-\pax^2 u^\epsilon$. We get
$$
\frac{1}{2}\frac{d}{dt}\|u^\epsilon\|_{\dot{H}^1}^2=I_1+I_2+I_3+I_4+I_5+I_6, 
$$
with
$$
I_1=\int_\TT\mu(u^\epsilon)\Lambda u^\epsilon\pax^2u^\epsilon dx,
\qquad
I_2=\int_\TT(\mu '(u^\epsilon))\pax u^\epsilon Hu^\epsilon\pax^2u^\epsilon dx,
$$
$$
I_3=-\int_\TT u^\epsilon(u^\epsilon-\langle u^\epsilon(t)\rangle)\pax^2u^\epsilon dx,
\qquad
I_4=-\int_{\TT}\pax u^\epsilon\pax v^\epsilon\pax^2u^\epsilon dx,
$$
$$
I_5=-\int_\TT ru^\epsilon(1-u^\epsilon)\pax^2 u^\epsilon dx,
\qquad
I_6=-\epsilon\|\pax^2u^\epsilon\|_{L^2}^2.
$$
As we have the bound \eqref{max}, using Assumption \ref{ass1} we get
$$
|\mu(u^\epsilon)|+ |\mu'(u^\epsilon)| \leq \mathcal{C}(s_0).
$$
Recall the following property of the Hilbert transform:
$$
\|H f\|_{L^p}\leq c_p\|f\|_{L^p},\;1<p<\infty.
$$
As a consequence, we find the bounds
$$
I_1\leq \mathcal{C} (s_0)\|u^\epsilon\|_{\dot{H}^1}\|\pax^2 u^\epsilon\|_{L^2},
$$
\[
I_2\leq \mathcal{C}(s_0)\|Hu^\epsilon\|_{L^4}\|\pax u^\epsilon\|_{L^4}\|\pax^2 u^\epsilon\|_{L^2}
\leq \mathcal{C}(s_0)c_4 \|u^\epsilon\|_{L^4}3\|u^{\epsilon}\|_{L^\infty}^{0.5}\|\pax^2 u^\epsilon\|_{L^2}^{1.5},
\]
$$
I_3\leq \|\pax^2 u^\epsilon\|_{L^2}\|u^\epsilon\|_{L^2}[\|u^\epsilon\|_{L^\infty}+\langle u^\epsilon(t) \rangle],
$$
$$
I_4\leq \frac{1}{2}\|u^\epsilon\|_{\dot{H}^1}^2[\|u^\epsilon\|_{L^\infty}+\langle u^\epsilon(t) \rangle],
$$
$$
I_5\leq r\|\pax^2 u^\epsilon\|_{L^2}\|u^\epsilon\|_{L^2}[\|u^\epsilon\|_{L^\infty}+1].
$$
Using Young's inequality and the $\epsilon$-dissipative term $I_6$, we have
\begin{equation}\label{H1odi}
\frac{d}{dt}\|u^\epsilon\|_{\dot{H}^1}^2\leq C_1(\epsilon)+C_2(\epsilon)\|u^\epsilon\|_{\dot{H}^1}^2, 
\end{equation}
where $C_1$ and $C_2$ depends also on $\mathcal{N}_1$, $s_0$, $r$ and $\gamma$. 

Let us assume now that for a given $\epsilon$ it holds  $T_\epsilon < T$. 
Using Gronwall's and Poincar\'e's inequalities in \eqref{H1odi}, we obtain existence of an approximate solution in $H^1$ up to  $T_\epsilon$. It is also bounded in view of \eqref{max}. Hence it is smooth by bootstrapping in  \eqref{eqa1approx}. Therefore it can be continued beyond  $T_\epsilon$, which consequently cannot be the maximal time of existence. We have $T_\epsilon \ge T$  and \eqref{L1}, \eqref{max} hold on $[0, T]$.

\subsection{Uniform estimates}
At this stage, we have  $\epsilon$-uniform estimates for $u^\epsilon(x,t)$ in $L^\infty_t L^p_x$ with $1\leq p\leq \infty$ on $[0, T]$ (see \eqref{max}).

Recalling \eqref{Gamma}, we have
$$
\int_\TT \pat u^\epsilon\log(u^\epsilon)dx=I_7+I_8+I_9+I_{10},
$$
with
$$
I_7=\int_\TT\pax u^\epsilon\gamma(u^\epsilon)Hu^\epsilon dx=\int_{\TT}\pax \Gamma(u^\epsilon)Hu^\epsilon dx=-\int_\TT \Gamma(u^\epsilon)\Lambda u^\epsilon dx,
$$
$$
I_8=-\int_\TT\pax u^\epsilon\pax v^\epsilon dx=\int_{\TT} u^\epsilon(u^\epsilon-\langle u^\epsilon\rangle) dx,
$$
$$
I_9=\int_\TT ru^\epsilon(1-u^\epsilon)\log(u^\epsilon)dx\leq 0,
$$
$$
I_{10}=-4\epsilon\int_{\TT}|\pax(\sqrt{u^\epsilon})|^2dx.
$$
Consequently, the evolution of the entropy functional \eqref{entropy} is given by
\begin{multline*}
\frac{d}{dt}\mathcal{F}(u^\epsilon(t))+\int_\TT \Gamma(u^\epsilon)\Lambda u^\epsilon dx+4\epsilon\int_{\TT}|\pax(\sqrt{u^\epsilon})|^2dx\\
=\|u^\epsilon-\langle u^\epsilon\rangle\|_{L^2}^2+\int_\TT ru^\epsilon(1-u^\epsilon)\log(u^\epsilon)dx.
\end{multline*}
Symmetrizing the integral $I_7$, we have
\begin{eqnarray*}
-I_7&=&\frac{1}{4\pi}\int_\TT\int_\TT\frac{u^\epsilon(x)-u^\epsilon(y)}{\sin^2\left(\frac{x-y}{2}\right)}\left(\Gamma(u^\epsilon(x))-\Gamma(u^\epsilon(y))\right)dxdy\\
&=&\frac{1}{4\pi}\int_\TT\int_\TT\int_0^1\frac{(u^\epsilon(x)-u^\epsilon(y))^2}{\sin^2\left(\frac{x-y}{2}\right)}\gamma(su^\epsilon(x)+(1-s)u^\epsilon(y))dsdxdy\\
&\geq&\frac{\delta}{4\pi}\int_\TT\int_\TT\frac{(u^\epsilon(x)-u^\epsilon(y))^2}{\sin^2\left(\frac{x-y}{2}\right)}dxdy = 
\delta\int_\TT\Lambda u^\epsilon u^\epsilon dx = \delta\|u^\epsilon\|_{\dot{H}^{0.5}}^2,
\end{eqnarray*}
where we have used \eqref{mu} of Assumption \ref{ass1}. Integrating in time, we get the uniform bound
$$
\mathcal{F}(u^\epsilon(t))+\delta\int_0^t\|u^\epsilon(s)\|_{\dot{H}^{0.5}}^2ds\leq C_3(\gamma,\mathcal{N}_1, s_0)t+\mathcal{F}(u^\epsilon_0),
$$
so
\begin{equation}\label{Hhalf}
\delta\int_0^t\|u^\epsilon(s)\|_{\dot{H}^{0.5}}^2ds\leq C_4(\gamma,\mathcal{N}_1, s_0)(t+1).
\end{equation}
For the derivative $\pat u^\epsilon$, using the duality pairing, we obtain
\begin{equation}\label{tder}
\|\pat u^\epsilon\|_{H^{-1.5}}\leq  \epsilon\|u^{\epsilon}\|_{H^{0.5}} + C (\gamma, s_0) 
\end{equation}
where we have used that $\|u^\epsilon\|_{L^2}$ controls $\|\pax v^\epsilon\|_{L^\infty}$. 

Let us sum up the obtained uniform bounds. Given  a finite, arbitrary $T\in (0, \infty)$ we have now the $\epsilon$-uniform bounds for
$$
u^\epsilon \text{ in } L^\infty(0,T; L^\infty(\TT))\cap L^2(0,T; H^{0.5}(\TT)),
$$ 
in view of \eqref{max} and \eqref{Hhalf} as well as for 
$$
\pat u^\epsilon  \text{ in } L^2(0,T;H^{-1.5}(\TT)).
$$ 
in view of \eqref{tder}, provided $\epsilon \le 1$.

For any $T < \infty$, applying the sequential $*$-weak compactness of spaces with a separable predual space, we have the existence of a limit function $u$ in appropriate spaces. Equipped with this $u$ we define $v$ using \eqref{eqa2}.

\subsection{Compactness} 
Applying Aubin-Lions's Lemma, we have the strong convergences (up to a subsequence) 
$$
\lim_{\epsilon\rightarrow 0}\int_0^T\|u^\epsilon(s)-u(s)\|_{L^2}^2ds=0.
$$
Notice that
$$
\lim_{\epsilon\rightarrow 0}\int_0^T|\langle u^\epsilon\rangle-\langle u\rangle|^2ds=0.
$$
Consequently, using
$$
\|\pax v^\epsilon-\pax v\|_{L^2}\leq \|u^\epsilon-u+\langle u\rangle-\langle u^\epsilon\rangle\|_{L^2}
$$
we get
\begin{equation}\label{vlimit}
 \lim_{\epsilon\rightarrow 0}\int_0^T\|\pax v^\epsilon-\pax v\|_{L^2}^2ds=0.
\end{equation}

\subsection{Passing to the limit} Here we show that \eqref{eqa1approx} is \eqref{eqa1} in the limit $\epsilon \to 0$, in the sense of  Definition \ref{defipe}.
We have in view of Assumption \ref{ass1} and \eqref{max}
$$
|\mu(u^\epsilon(x,t))-\mu(u(x,t))|\leq \mathcal{C}(\gamma,\mathcal{N}_\infty)|u^\epsilon(x,t)-u(x,t)|,
$$
and, consequently,
\begin{multline*}
\left|\int_0^T\int_{\TT}\pax \phi[\mu(u^\epsilon)-\mu(u)]H u^{\epsilon}dxds\right|\\
\leq \|\pax \phi\|_{L^\infty_tL^\infty_x}\|\mu(u^\epsilon)-\mu(u)\|_{L^2_tL^2_x}\|H u^\epsilon\|_{L^2_tL^2_x} \stackrel{\epsilon \to 0}{\longrightarrow}0.
\end{multline*}
Similarly we show
$$
\left|\int_0^T\int_{\TT}\pax \phi\mu(u)[H u^{\epsilon}-H u]dxds\right| \stackrel{\epsilon \to 0}{\longrightarrow} 0.
$$
Using \eqref{vlimit} we get
$$
\left|\int_0^T\int_\TT\pax\phi[u^\epsilon-u]\pax v^\epsilon dxds\right|\stackrel{\epsilon \to 0}{\longrightarrow}  0,\quad \left|\int_0^T\int_\TT\pax\phi u[\pax v^\epsilon-\pax v] dxds\right|  \stackrel{\epsilon \to 0}{\longrightarrow}  0.
$$
We can pass to the limit in the logistic term with the same ideas. Dealing with the term with laplacian is straightforward since it is linear. The only term left  is the one corresponding to the initial data. This term can be handled using the properties of mollifiers. 

We have proved that $(u,v)$ is a solution of \eqref{eqa1}-\eqref{eqa2} according to Definition \ref{defipe} and enjoys regularity properties from the statement of Theorem \ref{teo1}.

\section{Proof of Theorem \ref{teoRlarge}}\label{sec4}
Several steps are similar to those in the proof of Theorem \ref{teo1}, so we omit them. We consider the same approximate problems and we get the same $L^\infty_t L^1_x$ bound \eqref{L1}. Furthermore, these approximate problems have global solution in $H^1$. Consequently, we focus on the uniform estimates up to a fixed (but otherwise arbitrary) $0<T<\infty$.
\subsection{Pointwise bounds}
We use the same notation \eqref{notationmax} as before and we get that positivity is preserved. In particular, using
$$
-\Lambda u^\epsilon(\underline{x}_t,t)\geq -(u^\epsilon(\underline{x}_t,t)-\langle u^\epsilon(t)\rangle ),
$$
we can sharpen our bound. In this case, we have the ODE
\begin{eqnarray*}
\frac{d}{dt}\min_{x\in\TT} u^\epsilon(t)&=&\partial_t u^\epsilon(\underline{x}_t,t)\nonumber\\
&\geq&  (1-\gamma(u^\epsilon(\underline{x}_t))) u^\epsilon(\underline{x}_t)(u^\epsilon(\underline{x}_t)-\langle u^\epsilon(t)\rangle )+ru^{\epsilon}(\underline{x}_t)(1-u^{\epsilon}(\underline{x}_t))\nonumber\\
&\geq&  -\mu(u^\epsilon(\underline{x}_t)) (u^\epsilon(\underline{x}_t)-\langle u^\epsilon(t)\rangle )+u^{\epsilon}(\underline{x}_t)(1-\langle u^\epsilon(t)\rangle )\\
&\geq& -u^{\epsilon}(\underline{x}_t)\max\{1,\langle u_0\rangle\}, 
\end{eqnarray*}
so
$$
\min_{x\in\TT} u^\epsilon(t)\geq \text{ess\,inf}_{x\in\TT} u_0e^{-\max\{1,\langle u_0\rangle\} t}\geq 0\quad \forall \; 0\leq t\leq T<\infty.
$$
In particular, if the initial data is strictly positive,
\begin{equation}\label{strPos}
\min_{0\leq t\leq T}\min_{x\in\TT} u^\epsilon(t)\geq \text{ess\,inf}_{x\in\TT} u_0e^{-\max\{1,\langle u_0\rangle\} T} >0.
\end{equation}
We have to deal with the bound for $\|u^\epsilon(t)\|_{L^\infty}$. To this end we recall the inequality (see Lemma 1 in \cite{GO})
\begin{equation}\label{LuBO}
\Lambda u^\epsilon(\overline{x}_t) \geq \frac{u^\epsilon(\overline{x}_t)^2}{4\pi^2\langle u^\epsilon(t)\rangle}\geq \frac{u^\epsilon(\overline{x}_t)^2}{4\pi^2\max\{\langle u_0\rangle,1\}}.
\end{equation}
that is valid provided $u^\epsilon(\overline{x}_t) \ge 4 \langle u^\epsilon(t)\rangle $.
The ODI for $X (t) :=  u^\epsilon(\overline{x}_t)$ reads
\begin{eqnarray*}
\dot X
&\leq& -\mu(X)\Lambda X +X(X-\langle u^\epsilon(t)\rangle)+rX(1-X))\\
&\leq& X \left( \left(1 - r - \mu(X) \frac{ \Lambda X }{X^2} \right) X  +r \right)
\end{eqnarray*}
We proceed as before via a blowup alternative. 
Recall Assumption \ref{ass2}. In its context, choose $s_0\in\RR^+$ so large  that $\mu (s)\geq \delta$ for $s\ge s_0$ and that
\begin{equation}\label{Th2ch}
s_0 \ge \frac{2}{\pi}  \mathcal{N}_1, \quad s_0 \ge \frac{-2r}{1-r-\delta (4\pi^2\max\{\langle u_0\rangle,1\})^{-1}}. 
\end{equation}
$\mathcal{N}_1$ is given by \eqref{n1}. The second choice in \eqref{Th2ch} is possible thanks to $r +  \delta (4\pi^2\max\{\langle u_0\rangle,1\})^{-1} > 1$ assumed for our theorem.  We have the alternative
   \begin{itemize}
\item[(i)] Either $X (t)  \le s_0$ for all times.
\item[(ii)] Or there exists $t_0 > 0$ such  $X (t_0)  = s_0$ and $X$ crosses $s_0$ for the first time at $t_0$. 
\end{itemize}
It the latter case, we can use \eqref{LuBO} thanks to choice \eqref{Th2ch}  and \eqref{L1}. Hence we get from our ODI that $\dot X \le -rX
 < 0$, which excludes the case (ii).

\subsection{Uniform estimates}
We define
\begin{equation}\label{Mu}
\mathbb{M}(s):= \int_{0}^s\mu(y)dy.
\end{equation}
We test the equation \eqref{eqa1approx} against $u^\epsilon$ and and use \eqref{Mu}. We get
$$
\frac{1}{2}\frac{d}{dt} \int_\TT |u^\epsilon(t)|^2=\int_\TT \pax (\mathbb{M}(u^\epsilon) )  H u^\epsilon- \frac{1}{2}  \pax (|u^\epsilon|^2) \pax v^\epsilon+r |u^\epsilon(t)|^2 -r |u^\epsilon(t)|^3.
$$
After integration by parts and use of $v^\epsilon_{xx} = u^\epsilon - \langle u^\epsilon \rangle$ it yields
$$
\frac{1}{2}\frac{d}{dt} \int_\TT |u^\epsilon(t)|^2+ \int_\TT\mathbb{M}(u^\epsilon) \Lambda u^\epsilon + \left(r-\frac{1}{2}  \right) |u^\epsilon(t)|^3 =\left(r-\frac{\langle u^\epsilon(t)\rangle}{2}\right) \int_\TT  |u^\epsilon(t)|^2. 
$$
Hence
we obtain 
\[
\frac{d}{dt}\|u^\epsilon(t)\|_{L^2}^2+2\int_\TT\mathbb{M}(u^\epsilon) \Lambda u^\epsilon\leq\left(2r + \frac{\mathcal{N}_1}{2 \pi} \right)\|u^\epsilon(t)\|_{L^2}^2 + \left(\frac{1}{2}-r  \right) \|u^\epsilon(t)|\|_{L^3}^3,
\]
where we have used bound \eqref{L1} to control $\langle u^\epsilon(t)\rangle$. For $r \ge \frac{1}{2}$ the last term above provides extra dissipation, but for any $r \ge 0$ we can use our pointwise bounds to write
\begin{equation}\label{uniR}
\frac{d}{dt}\|u^\epsilon(t)\|_{L^2}^2+2\int_\TT\mathbb{M}(u^\epsilon) \Lambda u^\epsilon\leq   \frac{1}{2} s^3_0 + \left(2r + \frac{\mathcal{N}_1}{2 \pi} \right) s^2_0
\end{equation}
In order to extract $H^{0.5}$ information from \eqref{uniR} we symmetrize its second term analogously to $I_7$ of Section \ref{S2}, getting
\begin{equation*}
2\int_\TT\mathbb{M}(u^\epsilon) \Lambda u^\epsilon
=\frac{1}{2\pi}\int_\TT\int_\TT\int_0^1\frac{(u^\epsilon(x)-u^\epsilon(y))^2}{\sin^2\left(\frac{x-y}{2}\right)}\mu(su^\epsilon(x)+(1-s)u^\epsilon(y))dsdxdy\\
\end{equation*}
First, let's consider the case when $\mu$ verifies \eqref{mularger'}. Hence we have \eqref{strPos}. To simplify notation, we define
\begin{equation}\label{lower_rem}
0<s_1(T)=\text{ess\,inf}_{x\in\TT} u_0e^{-\max\{1,\langle u_0\rangle\} T}
\end{equation}
$$
\underline{\mu}_T=\min_{s_1\leq s\leq s_0} \mu(s),
$$
where $s_0$ is the $L^\infty$ bound from the previous subsection. We have
\begin{eqnarray*}
2\int_\TT\mathbb{M}(u^\epsilon) \Lambda u^\epsilon &\geq& \frac{\underline{\mu}_T}{2\pi} \int_0^1 \int_\TT\int_\TT\frac{(u^\epsilon(x)-u^\epsilon(y))^2}{\sin^2\left(\frac{x-y}{2}\right)} dxdy ds\\
&=&\frac{\underline{\mu}_T}{2\pi} \int_\TT\Lambda u^\epsilon u^\epsilon dx \\
&=&\frac{\underline{\mu}_T}{2\pi} \|u^\epsilon\|_{\dot{H}^{0.5}}^2
\end{eqnarray*}
Observe that $\underline{\mu}_T>0 $ as $\mu$ is positive for positive arguments.

If instead of \eqref{mularger'} we assume \eqref{mularger}, we have for $\delta >0$
\[
2\int_\TT\mathbb{M}(u^\epsilon) \Lambda u^\epsilon \ge
\frac{ \delta}{2\pi}  \|u^\epsilon\|_{\dot{H}^{0.5}}^2
\]
In any case we obtain via \eqref{uniR} the $\epsilon$-uniform bound for $u^\epsilon$ in $ L^2(0,T;H^{0.5}(\TT))$ for every fixed $0<T<\infty$.

All in all, for every fixed $0<T<\infty$, we have the uniform bounds
$$
u^\epsilon\in L^\infty(0,T;L^\infty(\TT))\cap L^2(0,T;H^{0.5}(\TT)).
$$
With these uniform estimates, we can follow along the lines of  the proof of Theorem \ref{teo1} and we conclude the result.
\begin{rem}\label{rem:stupido}
In fact, for Theorem \ref{teoRlarge} we need to impose assumptions on our semilinearity $\mu(s)$ only within the interval of existence of $u$. More precisely, this interval belongs to
\[I:= \left[\text{ess\,inf}_{x\in\TT} u_0e^{-\max\{1,\langle u_0\rangle\} T}, s_0 \right],
\]
where $s_0, \delta_0$ is a pair that satisfies
\[
s_0 \ge \frac{2}{\pi}  \mathcal{N}_1, \quad s_0 \ge \frac{-2r}{1-r-\delta_0 (4\pi^2\max\{\langle u_0\rangle,1\})^{-1}}, \quad \mu (s_0) > \delta_0 \ge 0. 
\]
For the lower end of $I$ see \eqref{lower_rem} and for the upper one -- \eqref{Th2ch}. In particular, $\mu$ can vanish on $I^c$. As the above condition on $s_0$ is implicit for $\delta_0 >0$, we have used in Theorem  \ref{teoRlarge} a more traceable assumption. For $\delta_0 =0$ we need $r>1$, but then 
\[I:= \left[\text{ess\,inf}_{x\in\TT} u_0e^{-\max\{1,\langle u_0\rangle\} T}, \frac{2r}{r-1} \right].
\]
and in fact $2$ above can be replaced with any $k>1$.
\end{rem}
\section{Proof of Corollary \ref{c1}}
For $r=0$ the conservation of mass in  \eqref{eqa1approx} gives
\[
\|u^{\epsilon}(t)\|_{L^1}\leq \|u_0\|_{L^1}.
\]
instead of \eqref{L1}. Consequently, \eqref{LuBO} reads now
\[
\Lambda u^\epsilon(\overline{x}_t) \geq \frac{u^\epsilon(\overline{x}_t)^2}{4\pi^2 \langle u_0 \rangle}.
\]
Hence, to follow the lines of proof of Theorem \ref{teoRlarge}, it suffices to assume \eqref{c:teoL'} instead of  \eqref{c:teoL}. This gives Corollary \ref{c1}.

\section{Proof of Theorem \ref{teo3}}\label{S3}
We consider the vanishing viscosity approximation of \eqref{eqa3} with $\gamma(x)\equiv 1$
$$
\pat u^\epsilon = -\pax(u^\epsilon Hu^\epsilon)+\epsilon\pax^2 u^{\epsilon},\;x\in\TT, t\in\RR^+.
$$
Notice that the solution to this equation verifies
$$
\min_x u^\epsilon(x,t)\geq \text{ess}\,\text{min}_x u_0>0,\;\max_x u^\epsilon(x,t) \leq \max_x u_0.
$$
By a direct computation, we have
$$
\frac{d}{dt}\mathcal{F}(u^\epsilon(t))+\mathcal{I}(u^\epsilon)+\epsilon\left\|\frac{\pax u^\epsilon}{\sqrt{u^\epsilon}}\right\|_{L^2}^2=0,
$$
where $\mathcal{F}$ is the entropy given by \eqref{entropy} and the Fisher's information $\mathcal{I}$ is
$$
\mathcal{I}(u^\epsilon)=\|\Lambda^{0.5}u^\epsilon\|_{L^2}^2,
$$
and verifies
\begin{multline*}
\frac{d}{dt}\mathcal{I}(u^\epsilon(t)) +2\epsilon\|\pax \Lambda^{0.5}u^\epsilon\|_{L^2}^2=-\|\sqrt{u^\epsilon}\Lambda u^\epsilon\|_{L^2}^2-\|\sqrt{u^\epsilon}\pax u^\epsilon\|_{L^2}^2\\
\leq -2\min_x u^\epsilon_0\mathcal{I}(u^\epsilon),
\end{multline*}
where the middle term follows from the Tricomi relation 
$$
H (H \pax u \pax u) = \frac{1}{2} ( ( H \pax u)^2 - (\pax u)^2),
$$ compare with \cite{GO}. In particular, by Gronwall inequality, we conclude that $u^\epsilon$ tends to the homogeneous state $\langle u^\epsilon_0 \rangle$ exponentially fast (recall that by assumption $\text{ess}\,\text{min}_x u_0 >0$). We also have
\begin{multline}\label{eqI}
\frac{d}{dt}\mathcal{F}(u^\epsilon(t))+\epsilon\left\|\frac{\pax u^\epsilon}{\sqrt{u^\epsilon}}\right\|_{L^2}^2=-\mathcal{I}(u^\epsilon(t))\\
\geq \frac{1}{2\min_x u^\epsilon_0}\left[\frac{d}{dt}\mathcal{I}(u^\epsilon(t))+2\epsilon\|\pax \Lambda^{0.5}u^\epsilon\|_{L^2}^2\right],
\end{multline}
and,
\begin{multline*}
\frac{d}{dt}\mathcal{F}(u^\epsilon(t))+\frac{\epsilon}{\min_x u^\epsilon_0}\|\pax u^\epsilon\|_{L^2}^2
\geq\frac{d}{dt}\mathcal{F}(u^\epsilon(t))+\epsilon\left\|\frac{\pax u^\epsilon}{\sqrt{u^\epsilon}}\right\|_{L^2}^2\\
\geq \frac{1}{2\min_x u^\epsilon_0}\frac{d}{dt}\mathcal{I}(u^\epsilon(t))+\frac{\epsilon}{\min_x u^\epsilon_0}\|\pax \Lambda^{0.5}u^\epsilon\|_{L^2}^2.
\end{multline*}
Due to Poincar\'e inequality, we get
$$
\frac{d}{dt}\mathcal{F}(u^\epsilon(t))\geq \frac{1}{2\min_x u^\epsilon_0}\frac{d}{dt}\mathcal{I}(u^\epsilon(t)).
$$
Equivalently,
$$
\int_t^\infty\frac{d}{dt}(-\mathcal{F}(u^\epsilon(t)))\leq \int_t^\infty\frac{-1}{2\min_x u^\epsilon_0}\frac{d}{dt}\mathcal{I}(u^\epsilon(t)),
$$
$$
-\mathcal{F}(u^\epsilon(\infty))+\mathcal{F}(u^\epsilon(t))\leq \frac{-1}{2\min_x u^\epsilon_0}\mathcal{I}(u^\epsilon(\infty))+\frac{1}{2\min_x u^\epsilon_0}\mathcal{I}(u^\epsilon(t)).
$$
As $\langle u^\epsilon_0\rangle=\langle u_0\rangle=1$, we have $\mathcal{F}(u^\epsilon(\infty))=0$, and we obtain
$$
\mathcal{F}(u^\epsilon(t))\leq \frac{1}{2\min_x u^\epsilon_0}\mathcal{I}(u^\epsilon(t)).
$$
Using \eqref{eqI},
$$
-\mathcal{F}(u^\epsilon(t))\geq \frac{1}{2\min_x u^\epsilon_0}\left[-\mathcal{I}(u^\epsilon(t))\right]=\frac{1}{2\min_x u^\epsilon_0}\left[\frac{d}{dt}\mathcal{F}(u^\epsilon(t))+\epsilon\left\|\frac{\pax u^\epsilon}{\sqrt{u^\epsilon}}\right\|_{L^2}^2\right]
$$
we conclude that
$$
-2\min_x u^\epsilon_0\mathcal{F}(u^\epsilon(t))\geq \frac{d}{dt}\mathcal{F}(u^\epsilon(t)),
$$
and
$$
\mathcal{F}(u^\epsilon(t))\leq \mathcal{F}(u^\epsilon_0)e^{-2(\text{ess}\,\text{min}_x u_0) \,t}.
$$
Hence, via $\mathcal{F}(u^\epsilon_0)\leq C(u_0)$ we have a uniform bound. Theorem 3.20 in \cite{Dacorogna} used for $a=b=0$ implies that the functional 
$$
f\mapsto \int_\TT f\log(f)-f+1dx
$$
is weakly lower semicontinuous in $L^2$. Hence
$$
\mathcal{F}(u(t))\leq \liminf_{\epsilon\rightarrow0} \mathcal{F}(u^\epsilon(t))\leq \mathcal{F}(u^\epsilon_0)e^{-2 (\text{ess}\,\text{min}_x u_0)\,t}.
$$
\section{Concluding remarks}
A preliminary computation (see the proofs of Theorem \ref{teo1} and \ref{teo2}) suggests that solutions to the case of subcritical diffusion $\Lambda^\alpha u$ with $\alpha \in (0, 1)$ should be bounded in $L^\infty_tL^\infty_x$, provided logistic dampening constant satisfies $r \ge 1$.

In the context of the critical diffusion  $\Lambda u$ and lack of the logistic term, the conjecture in \cite{bournaveas2010one} (see also \cite{AGM}) says that there should be a threshold mass that divides the global existence/blowup regimes. We are rather inclined against this hypothesis, along lack of the one-dimensional critical diffusion (\emph{i.e.} with threshold mass phenomenon) for the Smoluchowski-Poisson system with semilinear diffusion, compare \cite{CieslakLaurencot}. The authors show there that there is no one-dimensional analogue to the multidimensional phenomenon of the critical diffusion in the setting of semilinear, but not fractional diffusion. More precisely, for $d=1$ the system 

\begin{eqnarray*}
\pat u & = & \pax \left((1+u)^\frac{d-2}{d} \pax u +u\pax v \right), \\ 
\pax^2 v & = & u-\langle u \rangle,
\end{eqnarray*}

has bounded solutions for any initial mass.

\subsection*{Acknowledgments} JB is partially supported
by the National Science Centre (NCN) grant no. 2011/01/N/ST1/05411. RGB is partially supported by the grant MTM2011-26696 from the former Ministerio de Ciencia e Innovaci\'on (MICINN, Spain) and by the Department of Mathematics at University of California, Davis.

\bibliographystyle{abbrv}

\begin{thebibliography}{10}

\bibitem{AGM}
Y.~Ascasibar, R.~Granero-Belinch\'on, and J.~M. Moreno.
\newblock An approximate treatment of gravitational collapse.
\newblock {\em Physica D: Nonlinear Phenomena}, 262:71 -- 82, 2013.

\bibitem{BaeGranero}
H.~Bae and R. Granero-Belinch\'on
\newblock Global existence for some transport equations with nonlocal velocity.
\newblock {\em Submitted Arxiv preprint arXiv:1408.2768 [math.AP]}, 2014.


\bibitem{Bertozzi}
J.~Bedrossian, N.~Rodr{\'{\i}}guez, and A.~L. Bertozzi.
\newblock Local and global well-posedness for aggregation equations and
  {P}atlak-{K}eller-{S}egel models with degenerate diffusion.
\newblock {\em Nonlinearity}, 24(6):1683--1714, 2011.


\bibitem{Bi1}
P.~Biler.
\newblock {G}rowth and accretion of mass in an astrophysical model.
\newblock {\em Appl. Math.(Warsaw)}, 23(2):179--189, 1995.



\bibitem{Bi6}
P.~Biler, G.~Karch, P.~Lauren{\c{c}}ot, and T.~Nadzieja.
\newblock {T}he {$8\pi$}-problem for radially symmetric solutions of a
  chemotaxis model in the plane.
\newblock {\em {M}athematical {M}ethods in the {A}pplied {S}ciences},
  29(13):1563--1583, 2006.

\bibitem{BilNad94}
P.~Biler and T.~Nadzieja.
\newblock {E}xistence and nonexistence of solutions for a model of
  gravitational interaction of particles, i.
\newblock In {\em Colloq. Math}, volume~66, pages 319--334, 1994.




\bibitem{blanchet2010functional}
A.~Blanchet, E.~A. Carlen, and J.~A. Carrillo.
\newblock {F}unctional inequalities, thick tails and asymptotics for the
  critical mass {P}atlak-{K}eller-{S}egel model.
\newblock {\em {J}ournal of {F}unctional {A}nalysis}, 262(5):2142--2230, 2012.

\bibitem{blanchet2009critical}
A.~Blanchet, J.~Carrillo, and P.~Lauren{\c{c}}ot.
\newblock {C}ritical mass for a {P}atlak-{K}eller-{S}egel model with degenerate
  diffusion in higher dimensions.
\newblock {\em {C}alculus of {V}ariations and {P}artial {D}ifferential
  {E}quations}, 35(2):133--168, 2009.

\bibitem{BCM}
A.~Blanchet, J.~Carrillo, and N.~Masmoudi.
\newblock Infinite time aggregation for the critical {P}atlak-{K}eller-{S}egel
  model in $\mathbb{R}^2$.
\newblock {\em {C}ommunications on {P}ure and {A}pplied {M}athematics},
  61(10):1449--1481, 2008.

\bibitem{Dolbeault4}
A.~Blanchet, J.~Dolbeault, M.~Escobedo, and J.~Fern{\'a}ndez.
\newblock Asymptotic behaviour for small mass in the two-dimensional
  parabolic-elliptic {K}eller-{S}egel model.
\newblock {\em J. Math. Anal. Appl.}, 361(2):533--542, 2010.


\bibitem{BouCalChapter}
N.~Bournaveas and V.~Calvez.
\newblock  Kinetic models of chemotaxis. 
\newblock {\em Evolution equations of hyperbolic and Schršdinger type}, 41 -- 52, Progr. Math., 301, BirkhŠuser/Springer Basel AG, Basel, 2012.


\bibitem{bournaveas2010one}
N.~Bournaveas and V.~Calvez.
\newblock The one-dimensional {K}eller-{S}egel model with fractional diffusion
  of cells.
\newblock {\em {N}onlinearity}, 23(4):923, 2010.


\bibitem{BurczakCieslak}
J.~Burczak, T.~Cie{\'s}lak, and C.~Morales-Rodrigo.
\newblock Global existence vs. blowup in a fully parabolic quasilinear 1{D}
  {K}eller-{S}egel system.
\newblock {\em Nonlinear Anal.}, 75(13):5215--5228, 2012.

\bibitem{BurczakGranero}
J.~Burczak, and R.~Granero-Belinch\'on.
\newblock On a generalized Keller-Segel system in one spatial dimension.
\newblock {\em Submitted, Arxiv preprint arXiv:1407.2793 [math.AP]}, 2014.


\bibitem{Carrillo2}
V.~Calvez and J.~A. Carrillo.
\newblock Refined asymptotics for the subcritical {K}eller-{S}egel system and
  related functional inequalities.
\newblock {\em Proc. Amer. Math. Soc.}, 140(10):3515--3530, 2012.

\bibitem{Carrillo}
J.A. Carrillo, L.~F. Ferreira, and J. Precioso.
\newblock A mass-transportation approach to a one dimensional fluid mechanics
  model with nonlocal velocity.
\newblock {\em Adv. Math.}, 231(1):306--327, 2012.


\bibitem{CarrilloVazquez}
J.A. Carrillo, Y. Huang, M.C. Santos, and J.L. Vazquez.
\newblock Exponential Convergence Towards Stationary States for the 1D Porous Medium Equation with Fractional Pressure.
\newblock {\em arXiv:1407.4392 [math.AP]}, 2014.


\bibitem{CC}
A.~Castro and D.~C{\'o}rdoba.
\newblock {G}lobal existence, singularities and ill-posedness for a nonlocal
  flux.
\newblock {\em Advances in Mathematics}, 219(6):1916--1936, 2008.

\bibitem{CC2}
A. Castro and D. C{\'o}rdoba.
\newblock {S}elf-similar solutions for a transport equation with non-local
  flux.
\newblock {\em Chinese Annals of Mathematics, Series B}, 30(5):505--512, 2009.

\bibitem{CCCF}
D. Chae, A. C{\'o}rdoba, D. C{\'o}rdoba, and M.~A. Fontelos.
\newblock {F}inite time singularities in a {1D} model of the quasi-geostrophic
  equation.
\newblock {\em Advances in Mathematics}, 194(1):203--223, 2005.

\bibitem{CieslakLaurencot}
T. Cie\'slak and P. Lauren\,cot
\newblock Global existence vs. blowup in a one-dimensional Smoluchowski-Poisson system.
\newblock {\em Parabolic problems}, 95 -- 109, Progr. Nonlinear Differential Equations Appl., 80, BirkhŠuser/Springer Basel AG, Basel, 2011.

\bibitem{CieslakStinner}
T. Cie\'slak and C. Stinner
\newblock Finite-time blowup and global-in-time unbounded solutions to a parabolic-parabolic quasilinear Keller-Segel system in higher dimensions. 
\newblock {\em J. Differential Equations} 252, no. 10, 5832Ð5851, 2012.




\bibitem{Dacorogna} B. Dacorogna. 
Direct Methods in the Calculus of Variations. Springer, 2008. \newblock {\em}




%
%

\bibitem{Dolbeault2}
J.~Dolbeault and B.~Perthame.
\newblock Optimal critical mass in the two-dimensional {K}eller-{S}egel model
  in {$\Bbb R^2$}.
\newblock {\em C. R. Math. Acad. Sci. Paris}, 339(9):611--616, 2004.


\bibitem{escudero2006fractional}
C.~Escudero.
\newblock The fractional {K}eller-{S}egel model.
\newblock {\em Nonlinearity}, 19(12):2909, 2006.

%
%
\bibitem{GO}
R.~Granero-Belinch{\'o}n and R.~Orive-Illera.
\newblock An aggregation equation with a nonlocal flux.
\newblock {\em Nonlinear Analysis: Theory, Methods \& Applications}, 108(0):260
  -- 274, 2014.


\bibitem{Hillen3}
T.~Hillen and K.~J. Painter.
\newblock A user's guide to {PDE} models for chemotaxis.
\newblock {\em J. Math. Biol.}, 58(1-2):183--217, 2009.

\bibitem{Hillen1}
T.~Hillen, K.~J. Painter, and M.~Winkler.
\newblock Convergence of a cancer invasion model to a logistic chemotaxis
  model.
\newblock {\em Math. Models Methods Appl. Sci.}, 23(1):165--198, 2013.

%
\bibitem{jager1992explosions}
W.~J{\"a}ger and S.~Luckhaus.
\newblock {O}n explosions of solutions to a system of partial differential
  equations modelling chemotaxis.
\newblock {\em Trans. Amer. Math. Soc}, 329(2):819--824, 1992.

%

%
%




   
   

\bibitem{li2010exploding}
D.~Li, J.~Rodrigo, and X.~Zhang.
\newblock {E}xploding solutions for a nonlocal quadratic evolution problem.
\newblock {\em Revista Matematica Iberoamericana}, 26(1):295--332, 2010.


%

\bibitem{TaoWinkler}
Y. Tao and M.~Winkler. 
\newblock Boundedness in a quasilinear parabolic-parabolic Keller-Segel system with subcritical sensitivity. 
\newblock {\em J. Differential Equations 252 (2012), no. 1, 692 -- 715, 2012}


\bibitem{TelloWinkler}
J.~I. Tello and M.~Winkler.
\newblock A chemotaxis system with logistic source.
\newblock {\em Comm. Partial Differential Equations}, 32(4-6):849--877, 2007.

\bibitem{TelloWinkler2}
J.~I. Tello and M.~Winkler.
\newblock Stabilization in a two-species chemotaxis system with a logistic
  source.
\newblock {\em Nonlinearity}, 25(5):1413--1425, 2012.



\bibitem{WinklerLogistic}
M.~Winkler.
\newblock Boundedness in the higher-dimensional parabolic-parabolic chemotaxis
  system with logistic source.
\newblock {\em Comm. Partial Differential Equations}, 35(8):1516--1537, 2010.





\end{thebibliography}

\end{document}